\pdfoutput=1
\documentclass[11pt]{article}
\usepackage{graphicx,amsmath,amscd,amssymb,amsthm}
\usepackage{fourier}
\usepackage{listings}
\usepackage[margin=1in]{geometry}
\usepackage[sort&compress, comma, square, numbers]{natbib}
\usepackage[colorlinks,linkcolor=black,bookmarksopen,
bookmarksnumbered,citecolor=black,urlcolor=black]{hyperref}
\usepackage[small,bf]{caption}
\usepackage{enumerate}
\usepackage{color}
\usepackage{titlesec}
\usepackage{multirow}
\usepackage{array}
\usepackage{authblk}

\graphicspath{{figs/darcy/}}
\DeclareMathOperator{\grad}{grad}

\def\div{\operatorname{div}}

\def\d{\operatorname{d}}

\DeclareMathOperator{\dd}{d} 
\DeclareMathOperator{\hodge}{\ast}

\DeclareMathOperator{\boundary}{\partial}

\newtheorem*{theorem*}{Theorem}

\newtheorem*{proposition*}{Proposition}

\newtheorem*{lemma*}{Lemma}

\newtheorem*{claim*}{Claim}

\newtheorem*{axiom*}{Axiom}

\newtheorem*{conjecture*}{Conjecture}

\newtheorem*{corollary*}{Corollary}

\theoremstyle{definition}

\newtheorem*{definition*}{Definition}

\newtheorem*{example*}{Example}

\newtheorem*{exercise*}{Exercise}
\newtheorem*{recall*}{Recall}

\theoremstyle{remark}

\newtheorem*{note*}{Note}

\newtheorem*{remark*}{Remark}

\newtheorem*{notation*}{Notation}

\newtheorem*{question*}{Question}

\newtheorem*{fact*}{Fact}

\theoremstyle{definition}

\theoremstyle{remark}

\allowdisplaybreaks

\begin{document}

\title{\vspace*{-1.4cm} Numerical Experiments for Darcy Flow on a
  Surface \\Using Mixed Exterior Calculus Methods}

\author{ Anil N.~Hirani\thanks{email:
    \href{mailto:hirani@cs.illinois.edu}{hirani@cs.illinois.edu};
    URL: \url{http://www.cs.illinois.edu/hirani}}\quad}
\author{\quad Kaushik~Kalyanaraman}
\affil{Department of Computer
  Science, University of Illinois at Urbana-Champaign}

\date{}
\maketitle

\begin{abstract}
  There are very few results on mixed finite element methods on
  surfaces. A theory for the study of such methods was given recently
  by Holst and Stern, using a variational crimes framework in the
  context of finite element exterior calculus. However, we are not
  aware of any numerical experiments where mixed finite elements
  derived from discretizations of exterior calculus are used for a
  surface domain. This short note shows results of our preliminary
  experiments using mixed methods for Darcy flow (hence scalar
  Poisson's equation in mixed form) on surfaces. We demonstrate two
  numerical methods. One is derived from the primal-dual Discrete
  Exterior Calculus and the other from lowest order finite element
  exterior calculus. The programming was done in the language Python,
  using the PyDEC package which makes the code very short and easy to
  read. The qualitative convergence studies seem to be promising.
\end{abstract}

\section{Introduction}
In this short note we present results from some preliminary numerical
experiments for Darcy flow (hence scalar Poisson's equation in mixed
form) on a surface. We use two mixed methods. One is derived from
primal-dual Discrete Exterior Calculus (DEC) and adapted
from~\cite{HiNaCh2008}. 
The other method uses the $\mathcal{P}_1^- \Lambda^1$ elements,
i.e., Whitney 1-forms for the fluxes, and piecewise constants for the
pressures. We'll refer to this second method as being derived from
finite element exterior calculus (FEEC)~\cite{ArFaWi2006, ArFaWi2010}
or call it the Whitney solution.
Switching between the DEC and FEEC methods requires a 1 line change in
our code, which is written in Python using PyDEC \cite{BeHi2011}. The
only difference between the two methods in our code is that the former
method uses a primal-dual DEC Hodge star, and the latter uses Whitney
Hodge star~\cite{BeHi2011}.  The surface used is a hemisphere with a
hole punched out around north pole, and we will refer to this surface
as an annular hemisphere. The boundary conditions are Neumann for
pressure, hence specified as flux on the boundary. The pressure is
fixed arbitrarily at a point to make the problem well-posed. The use
of the Whitney Hodge star is equivalent to a solution using finite
element exterior calculus. Recall that the equations of Darcy flow on
domain $M$ are:
\begin{alignat*}{2}
v +  \grad p & = 0 && \quad \text{on } M \\
\div v & = 0 &&  \quad \text{on } M \\
v \cdot \hat{n} & = \psi && \quad \text{on } \partial M
\end{alignat*}
where we have assumed the ratio of permeability and viscosity to be 1
and assumed the absence of sources and sinks in the domain. The
unknowns are the pressure $p$ and velocity $v$, in $M$. In mixed FEM
literature, the preferred variable names are $u$ and $\sigma$,
respectively.

\section{Exterior Calculus Mixed Methods}
In exterior calculus notation, taking $p$ to be a 0-form the Darcy
flow equations can be written in terms of the 1-form proxy for the
vector field, and the exterior derivative of the pressure. The
divergence is replaced by the codifferential. In particular, on domain
$M$, the Darcy equations from above become
\begin{align*}
v^{\flat} + \d p &= 0 \, , \\
-\delta v^{\flat} &= 0 \, ,
\end{align*}
where $\flat$ is the isomorphism between vector fields and 1-forms and
$\delta$ is the codifferential~\cite{AbMaRa1988}. 
Applying a Hodge star to the first equation and using $\delta
= -\hodge \d \hodge v^{\flat}$ we have
\begin{align*}
\hodge v^{\flat} + \hodge \d p &= 0 \, , \\
\d \hodge v^{\flat} &= 0 \, .
\end{align*}
We will set $\sigma = \hodge v^{\flat}$ which is a flux, and use that
as an unknown. The above equations for a 2-dimensional domain are
equivalent to
\begin{align*}
\hodge \sigma - \d p &= 0 \, , \\
\d \sigma &= 0 \, .
\end{align*}
Following~\cite{HiNaCh2008} we discretize this
as
\[
\begin{bmatrix}
-\hodge_1 & \dd_1^T \\
\phantom{-}\dd_1 & 0
\end{bmatrix}%
\begin{bmatrix}
\sigma \\ 
p
\end{bmatrix} = %
\begin{bmatrix}
0 \\
0
\end{bmatrix}\, ,
\]
where the coboundary matrix $\dd_1$ is the transpose of the boundary
matrix $\boundary_2$ from triangles to edges. This linear system is
made nonsingular by fixing pressure at a point and adjusting the
system for this. For more details on the above derivation
see~\cite{HiNaCh2008}.  When $\hodge_1$ is the mass matrix for Whitney
1-forms, i.e., Whitney Hodge star, one gets a finite element exterior
calculus type Whitney solution. When $\hodge_1$ is the primal-dual DEC
Hodge star, one gets the DEC method of \cite{HiNaCh2008}. Both these
methods are used in this paper.

In the code, a \texttt{simplicial\_complex} object of PyDEC is created
after reading the mesh from files. If \texttt{sc} is the name of this
object, then the choice between Whitney or DEC solution requires only
specifying the appropriate Hodge star matrix in the Python code. For
the Whitney Hodge star, which yields the Whitney solution, one writes
\begin{lstlisting}[firstline=61,lastline=61]
    star1 = whitney_innerproduct(sc, 1)
\end{lstlisting} and for
the primal-dual DEC Hodge star, which yields the DEC solution, one
writes 
\begin{lstlisting}[firstline=59, lastline=59]
    star1 = sc[1].star
\end{lstlisting}
The matrix of the system above is assembled by the code
\begin{lstlisting}[firstline=62,lastline=63]
A = bmat([[(-mu/k)*star1, sc[1].d.T],
          [sc[1].d, None]], format='csr')
\end{lstlisting} where
\texttt{csr} refers to the compressed sparse row format, one of many
sparse matrix formats available via SciPy. The variables \texttt{mu}
and \texttt{k} refer to the viscosity and permeability and for this
note we are assuming this ratio to be 1. The main programming effort
in the rest of the code is in the quadratures needed for the exact
flux and in setting the boundary conditions.

Thus the numerical methods use flux as a 1-cochain associated with the
edges, and pressure associated with the triangles. In DEC the pressure
is a dual 0-cochain associated with the circumcenters of the triangles
of the mesh which needs to be Delaunay. In the Whitney solution the
pressure is considered constant on each triangle. For visualization,
we interpolate the flux using a Whitney map and sample it at the
barycenters of the triangles and rotate it by $\pi/2$.

\section{Planar Annulus as a Preparatory Exercise}
To motivate the solution on the annular hemisphere we first discuss
the case of planar annulus.  Consider a planar annulus with origin as
center, inner radius $r_0$ and outer radius $r_1$. The fluid enters
from the inner boundary and exits through the outer boundary. We will
take the direction of the flow to be normal to the boundary and speed
to be constant along the inner boundary. At the inflow boundary
velocity vectors are pointing into the domain, and at the outflow
boundary they are pointing out of the domain. By symmetry, the
velocity at any point will be directed radially pointing away from the
origin. The speed along a circle of radius $r$ will be constant. The
magnitude of velocity can be determined by using divergence theorem as
follows.

Let $D_r$ be an open disk of radius $r \, , r_0 < r \leq r_1$, and let
$M_r = M \cap D_r$ be the annulus with the inner radius $r_0$ and
outer radius $r$. Since there are no sources and sinks (i.e., $\div v
= 0 \text{ in } M$), we have
\[
0 = \int\limits_{M_r} \div v \, dx = \int\limits_{\partial M_r} v \cdot
\hat{n} \, d\gamma \, ,
\]
where $d\gamma$ is the measure on the boundary. If $\Gamma_0$ is
the inner circle and $\Gamma$ the outer circle, we have
\begin{equation} \label{eq:divbndry}
\int\limits_{\Gamma_0} v(x, y) \cdot \hat{n}(x, y) \, d\gamma =
\int\limits_{\Gamma}  v(x, y) \cdot \hat{n}(x, y) \, d\gamma \, ,
\end{equation}
Since $v \cdot \hat{n}$ is a constant on any circle, let us call the
first integrand $S(r_0)$ and the second one $S(r)$. Then we have
$S(r_0) \, 2 \pi r_0 = S(r) \, 2 \pi r$, or 
\[
S(r) = S(r_0) \dfrac{r_0}{r} \, .
\]

The norm of the velocity on outer circle of the original annulus is
$S(r_0) r_0/r_1$. From the expression for speed we can determine the
pressure at any $r$ as follows. Along a ray from origin, we have a one
dimensional problem and
\[
S(r) = -\dfrac{dp}{dr} \, .
\]
Thus $p(r)$ along any ray is given by
\[
p(r) = -\int\limits_{r_0}^{r} S(\xi) \, d\xi + C_0 = -S(r_0) r_0 
\int\limits_{r_0}^r\dfrac{1}{\xi} \, d\xi + C_0 \, ,
\]
where $C_0$ is the arbitrary pressure at $r_0$. Thus
\[
p(r) = -S(r_0) r_0 \left. \ln r \right \vert_{r_0}^r + C_0 = -S(r_0)
r_0 \ln r + C
\]
where $C = C_0 + S(r_0) r_0 \ln r_0$, hence just another arbitrary
constant. In our code we took $C_0 = 0$, and hence $C = S(r_0) r_0 \ln
r_0$. In fact since the experiments use an annulus with $r_0=1$ even
$C = 0$ in these experiments.

\subsection*{Numerical results for planar annulus}
\begin{figure}[hb]
  \centering
  \includegraphics[scale=0.8, trim=0in 0in 0in 0in, clip]
  {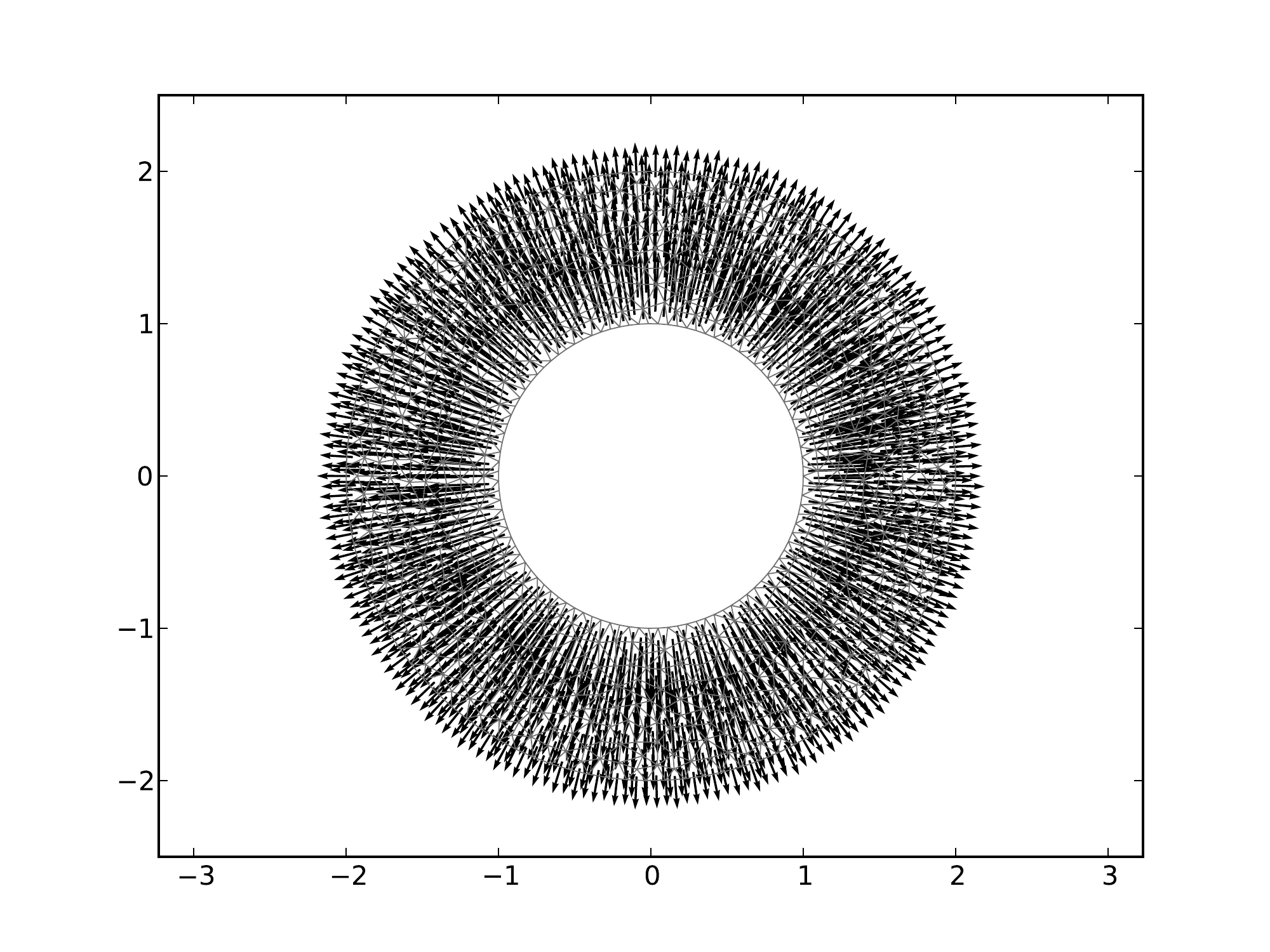}
  \caption{The velocity field corresponding to the flux on a planar
    annulus computed as a Whitney solution. The boundary condition is
    Neumann and pressure is fixed arbitrarily at a single point. The
    inflow is constant speed entry at inner boundary and constant
    speed exit at outer boundary whose value can be computed using the
    divergence theorem. The profiles for speed and pressure are in
    Figure~\ref{fig:annlsFEEC}. The planar annulus here is used as a
    warmup exercise for the surface experiments on an annular
    hemisphere.}
  \label{fig:whtnyvctrfld}
\end{figure}

\begin{figure}[ht]
  \begin{tabular}{m{0.475\textwidth}m{0.475\textwidth}}
    \hspace{-0.04\textwidth}\includegraphics[scale=0.45]
    {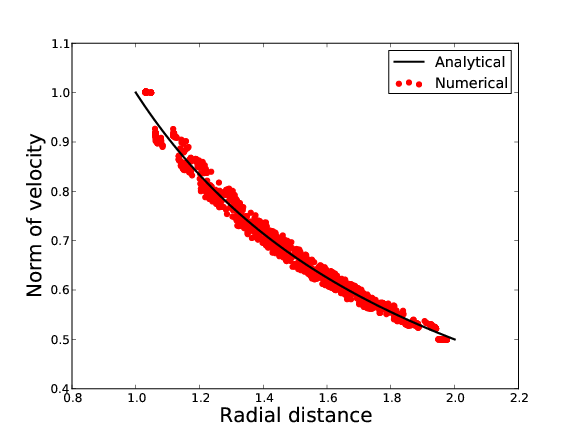} & \hspace{-0.02\textwidth}
    \includegraphics[scale=0.45]
    {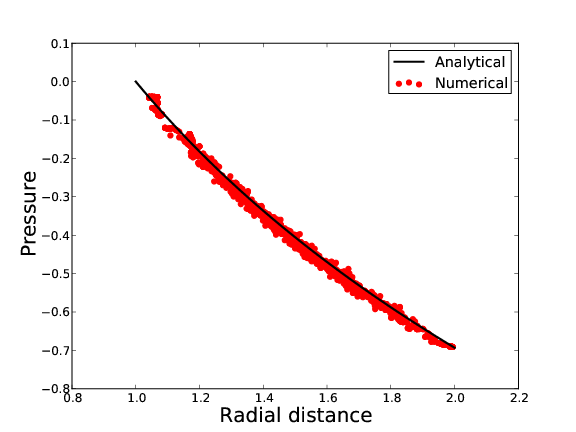}
  \end{tabular}
  \caption{Numerical solution for planar annulus computed as a Whitney
    solution, compared with analytical solution. Left plot shows a
    scatter plot of the speeds compared with the analytical values
    over a ray from the inner to outer boundary of the annulus. Right
    plot is a similar plot for the pressures.}
  \label{fig:annlsFEEC}
\end{figure}

\begin{figure}[hb]
  \begin{tabular}{m{0.475\textwidth}m{0.475\textwidth}}
    \hspace{-0.04\textwidth} 
    \includegraphics[scale=0.45]{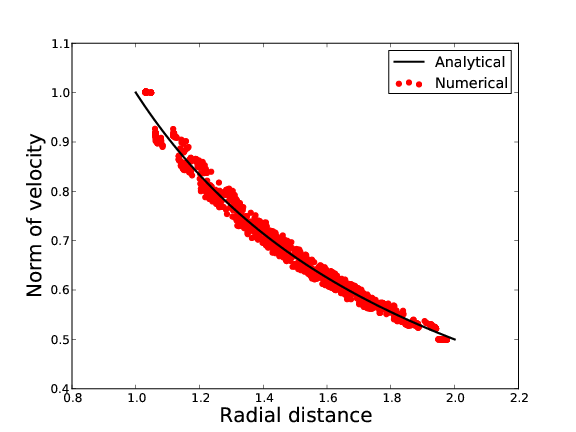} & 
    \hspace{-0.02\textwidth}
    \includegraphics[scale=0.45]{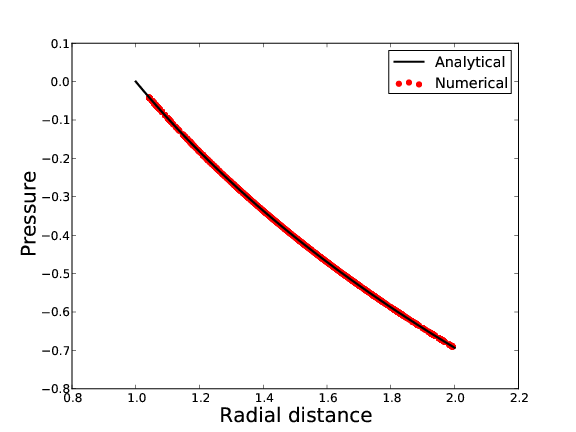}
  \end{tabular}
  \caption{Numerical results for planar annulus computed using
    DEC. The two plots are similar to the ones in
    Figure~\ref{fig:annlsFEEC}.}
  \label{fig:annlsDEC}
\end{figure}

\begin{figure}[p]
  \vspace{-0.5in}
  \begin{tabular}{m{0.475\textwidth}m{0.475\textwidth}}
    \hspace{-0.04\textwidth} \includegraphics[scale=0.45]
    {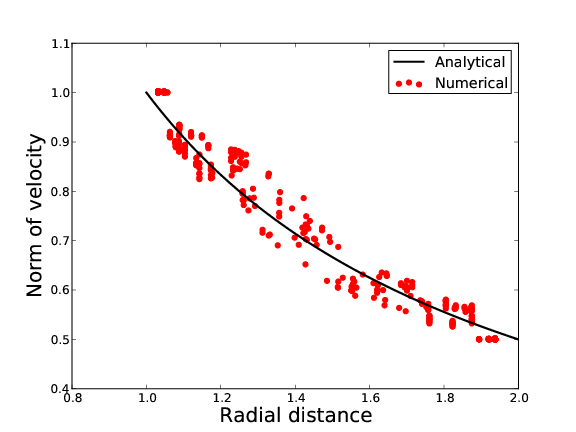}  & \hspace{-0.02\textwidth}
    \includegraphics[scale=0.45]{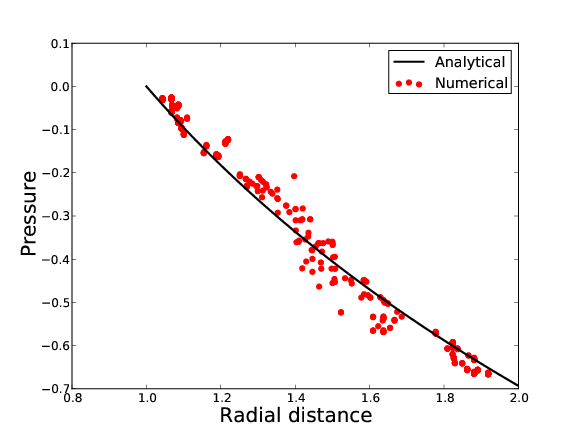}
    \\
    \hspace{-0.04\textwidth} \includegraphics[scale=0.45]
    {annulus/1516/whitney/speeds.png}  & \hspace{-0.02\textwidth}
    \includegraphics[scale=0.45]{annulus/1516/whitney/pressures.png}
    \\
    \hspace{-0.04\textwidth} \includegraphics[scale=0.45]
    {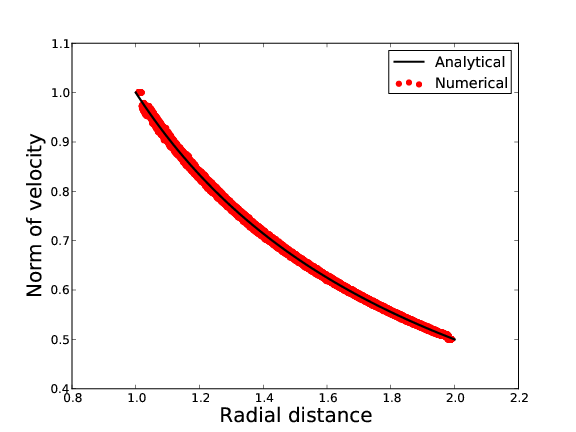}  & \hspace{-0.02\textwidth}
    \includegraphics[scale=0.45]{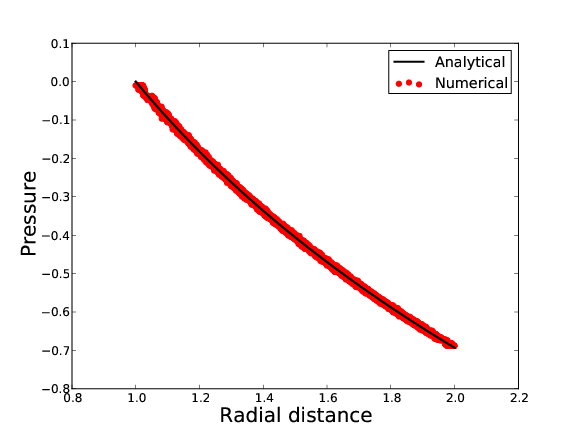}
  \end{tabular}
  \caption{A qualitative indication of convergence of the Whitney
    solution for the planar annulus on a series of Delaunay meshes
    similar to the one shown in Figure~\ref{fig:whtnyvctrfld}. Similar
    to Figure~\ref{fig:annlsFEEC}, the left panels are speed plots and
    the right ones are pressure plots. The rows of figures from top to
    bottom are profiles on increasingly finer meshes with 484, 1516
    and 14685 elements, respectively.}
  \label{fig:annlsFEECcnvrgnc}
\end{figure}

\begin{figure}[p]
  \vspace{-0.5in}
  \begin{tabular}{m{0.475\textwidth}m{0.475\textwidth}}
    \hspace{-0.04\textwidth} \includegraphics[scale=0.45]
    {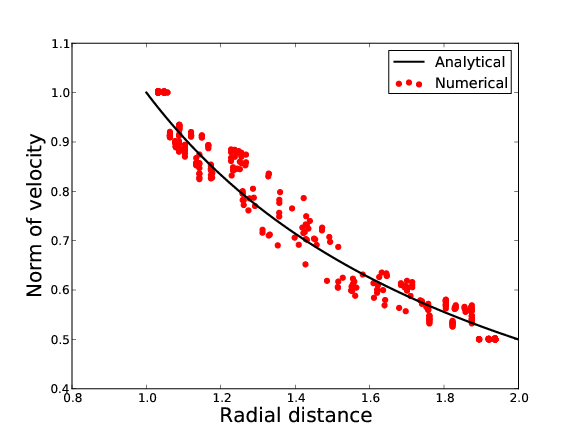}  & \hspace{-0.02\textwidth}
    \includegraphics[scale=0.45]{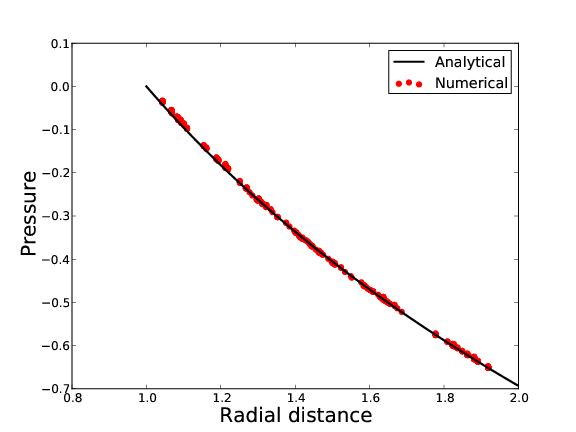}
    \\
    \hspace{-0.04\textwidth} \includegraphics[scale=0.45]
    {annulus/1516/dec/speeds.png}  & \hspace{-0.02\textwidth}
    \includegraphics[scale=0.45]{annulus/1516/dec/pressures.png}
    \\
    \hspace{-0.04\textwidth} \includegraphics[scale=0.45]
    {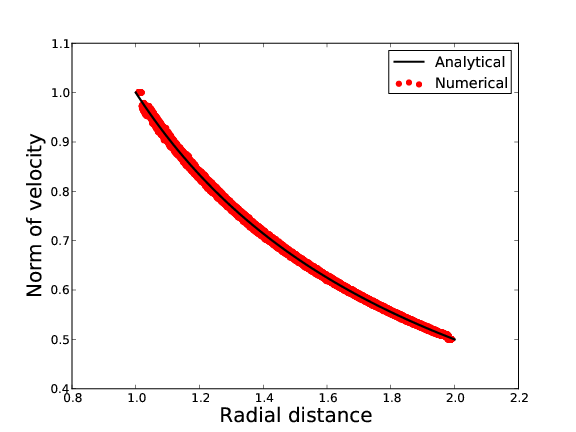}  & \hspace{-0.02\textwidth}
    \includegraphics[scale=0.45]{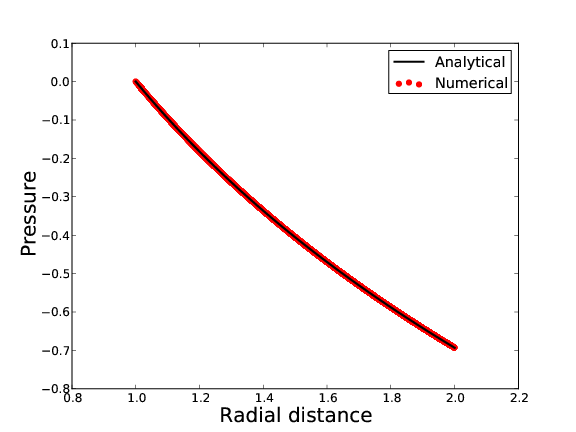}
  \end{tabular}
  \caption{A qualitative indication of convergence of DEC solution for
    the planar annulus on a series of Delaunay meshes similar to the
    one shown in Figure~\ref{fig:whtnyvctrfld}. Similar to
    Figure~\ref{fig:annlsFEEC}, the left panels are speed
    plots and the right ones are pressure plots. The rows of figures
    from top to bottom are profiles on increasingly finer meshes with
    484, 1516 and 14685 elements, respectively.}
  \label{fig:annlsDECcnvrgnc}
\end{figure}

We computed the pressure as a 0-cochain and flux as a 1-cochain. In
the DEC solution the pressure is a dual 0-cochain associated with the
circumcenters even if the circumcenter is outside the triangle, as
long as it is Delaunay~\cite{BeHi2011}. In the Whitney solution
solution, pressure is assumed constant on each triangle.  The results
of the annulus experiments are shown for the Whitney solution in
Figures~\ref{fig:whtnyvctrfld} and~\ref{fig:annlsFEEC} and for DEC in
Figure~\ref{fig:annlsDEC}.  Figure~\ref{fig:whtnyvctrfld} shows the
velocity vector field sampled at the barycenter of each triangle.
Figures~\ref{fig:annlsFEEC} and~\ref{fig:annlsDEC} show the speed and
pressure as a function of $r$ for the two different methods. In both,
the Whitney solution plots and the DEC plots, for the purpose of
visualization, the location for the pressure in a triangle was
associated with its circumcenter. That is the value of $r$ used for a
particular computed pressure was that of the circumcenter of the
triangle. The flux was interpolated using a Whitney map and sampled at
the barycenter and thus the $r$ value of the barycenter was used for
drawing the speed plots. Qualitative studies of convergence are in
Figure~\ref{fig:annlsFEECcnvrgnc} and~\ref{fig:annlsDECcnvrgnc}. We
started with a coarse mesh of 484 triangles and refined it by
subdividing each triangle into 4 congruent triangles to obtain meshes
with 1,516 and 14,685 triangles.

\section{Annular Hemisphere}
\begin{figure}[p]
  \vspace{-0.5in}
  \centering
  \begin{tabular}{c}
    \includegraphics[scale=0.7, trim=3in 1in 3in 2in, clip]
    {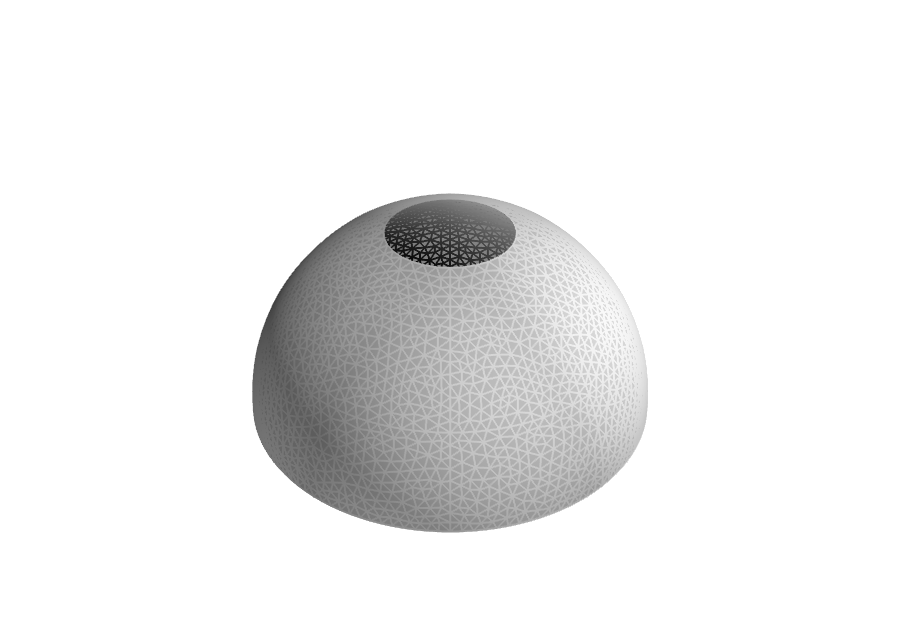}
    \\
    \includegraphics[scale=0.7, trim=3in 1in 3in 2in, clip]
    {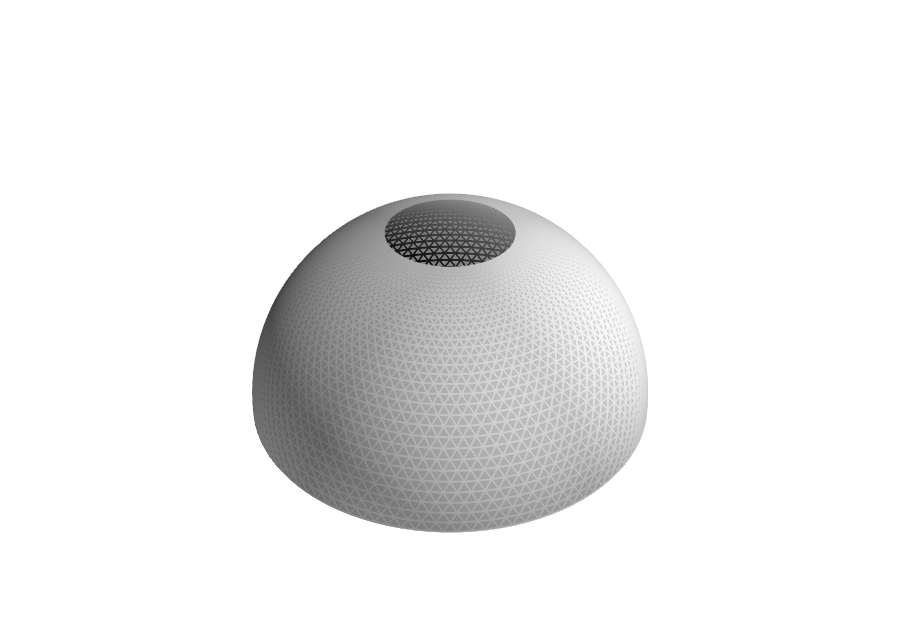}
  \end{tabular}
  \caption{Representative triangulations of the annular hemisphere.
    The top mesh is generated from the mesh generation software gmsh,
    and the bottom mesh is a well-centered
    triangulation~\cite{VaHiGuRa2010a} obtained using Python code. The
    mesh in the top figure consists of 4928 triangles. The speed and
    pressure profiles for Whitney solution obtained using a similar
    mesh consisting of 138536 triangles is shown in
    Figure~\ref{fig:hmsphrhlFEEC}. The bottom mesh consists of 6600
    triangles. The speed and pressure profiles for DEC solution
    obtained using a similar mesh consisting of 135900 triangles is
    shown in Figure~\ref{fig:hmsphrhlDEC}.}
  \label{fig:hmsphrhlmsh}
\end{figure}

Now we will consider a hemispherical domain with a circular hole
punched out around north pole. The bottom boundary sits on the
$xy$-plane and the $z$-axis points up. Two of the triangulations of
the domain that we use are shown in Figure~\ref{fig:hmsphrhlmsh}. The
inflow is now from the top boundary and the outflow at the
equator. Both are tangential to the surface and normal to the
boundary.

The derivation of the analytical solution closely follows the
derivation for the planar annulus. By symmetry the speed at a fixed
latitude will be constant, with velocity normal to the latitude and
tangential to the sphere. The divergence theorem can be used as in the
planar case and the shape of the surface between the two latitudes is
irrelevant. Consider any point on the surface where the velocity and
pressure are desired. Let $r$ be the horizontal distance of this point
from the $z$-axis. The derivation for speed follows the one of planar
annulus, with $\Gamma_0$ being the top boundary and $\Gamma_r$ being
the latitude passing through the point that is distance $r$ from the
$z$-axis.

For deriving the pressure solution, we will use spherical coordinates
with the convention that $\theta$ is measured from $z$-axis and $\phi$
is measured from the $x$-axis. Since the radius of the hemisphere is 1,
we will use $r$ to stand for the distance from $z$-axis. By symmetry,
the only nonzero component for $v$ is in the $\theta$-direction. We
will call this component $S(\theta)$. By the planar annulus
derivation, 
\[
S(\theta) = S(\theta_0) \, \dfrac{\sin \theta_0}{\sin \theta} \, ,
\]
where $\theta_0$ corresponds to the latitude at $r_0$ distance from
the $z$-axis. Thus,
\[
\frac{d p}{d \theta} = -S(\theta_0) \, \dfrac{\sin \theta_0}{\sin
\theta} \, ,
\]
which implies 
\[
p(\theta) = -S(\theta_0) \, \sin \theta_0 \int\limits_{\theta_0}^{\theta}
\dfrac{1}{\sin \xi} d \xi + C_0 \, ,
\]
for an arbitrary constant $C_0$. Thus,
\[
p(\theta) = S(\theta_0) \, \sin \theta_0 \; \ln \left( \dfrac{1 + \cos
    \theta}{\sin \theta} \right) + C \, ,
\]
and we use
\[
C = -S(\theta_0) \, \sin \theta_0 \; \ln  \left( \dfrac{1 + \cos
    \theta_0}{\sin \theta_0} \right) + C_0 \, ,
\]
with $C_0$ set to 0 in our code.

\subsection*{Numerical results for annular hemisphere}

\begin{figure}[!htp]
  \begin{tabular}{m{0.475\textwidth}m{0.475\textwidth}}
    \hspace{-0.04\textwidth} \includegraphics[scale=0.45]
    {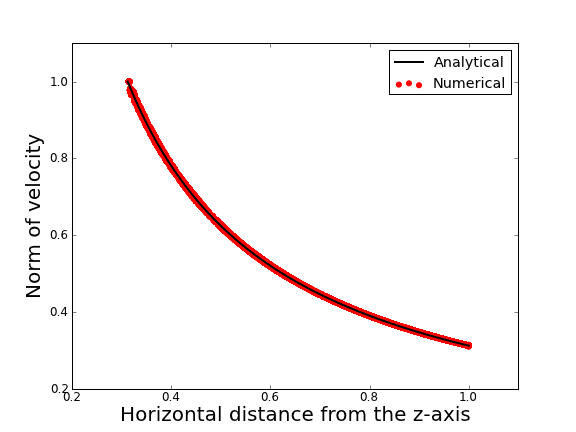}  & \hspace{-0.02\textwidth}
    \includegraphics[scale=0.45]{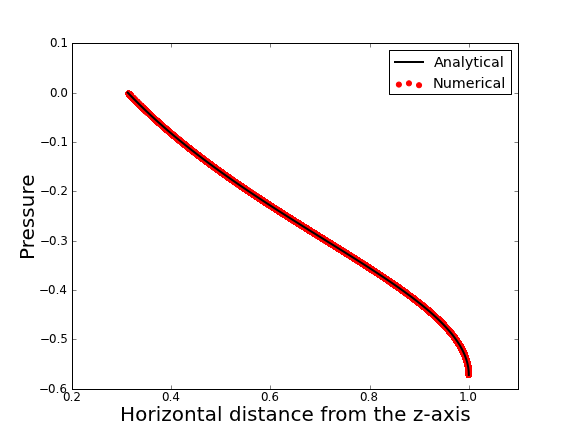}
  \end{tabular}
  \caption{Whitney solution for an annular hemisphere of the type
    shown in top of Figure~\ref{fig:hmsphrhlmsh}. Left figure here
    shows a scatter plot of the speeds compared with the analytical
    values over a longitude from the top to bottom boundary. Right
    figure shows a similar plot for the pressures.}
  \label{fig:hmsphrhlFEEC}
\end{figure}

\begin{figure}[!htp]
  \begin{tabular}{m{0.475\textwidth}m{0.475\textwidth}}
    \hspace{-0.04\textwidth} \includegraphics[scale=0.45]
    {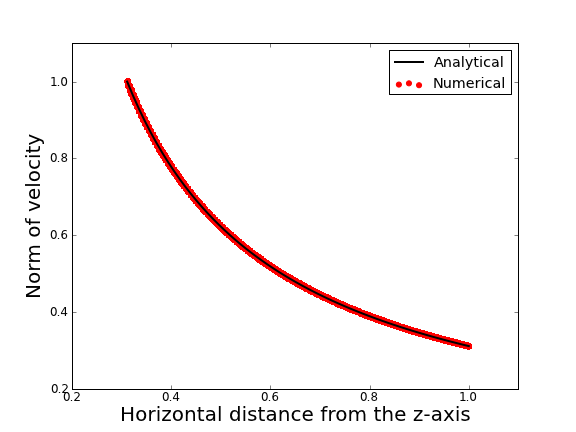}  & \hspace{-0.02\textwidth}
    \includegraphics[scale=0.45]{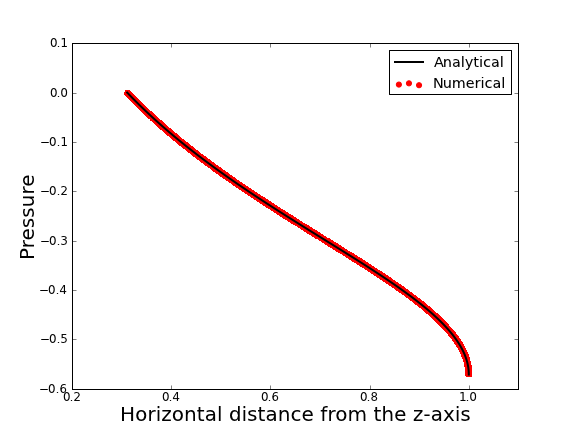}
  \end{tabular}
  \caption{DEC solution for a well-centered triangulation of a
    punctured hemisphere of the type shown in bottom of 
    Figure~\ref{fig:hmsphrhlmsh}. The two plots are similar to the
    ones in Figure~\ref{fig:hmsphrhlFEEC}.}
  \label{fig:hmsphrhlDEC}
\end{figure}

\begin{figure}[p]
 \begin{tabular}{m{0.475\textwidth}m{0.475\textwidth}}
    \hspace{-0.04\textwidth} 
    \includegraphics[scale=0.45]{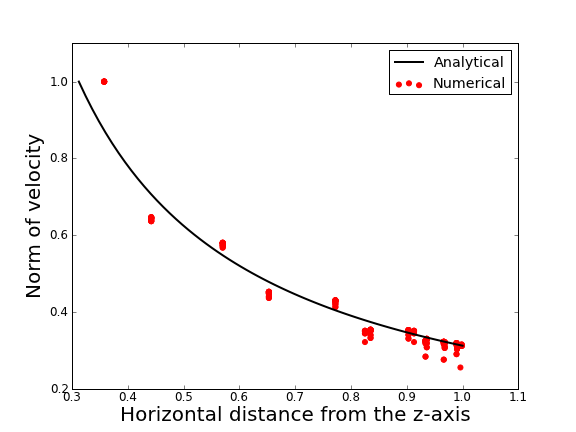}  & 
    \hspace{-0.02\textwidth}
    \includegraphics[scale=0.45]{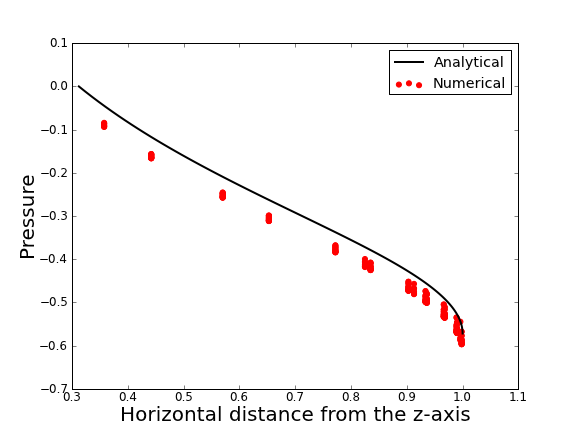} \\
    \hspace{-0.04\textwidth} 
    \includegraphics[scale=0.45]{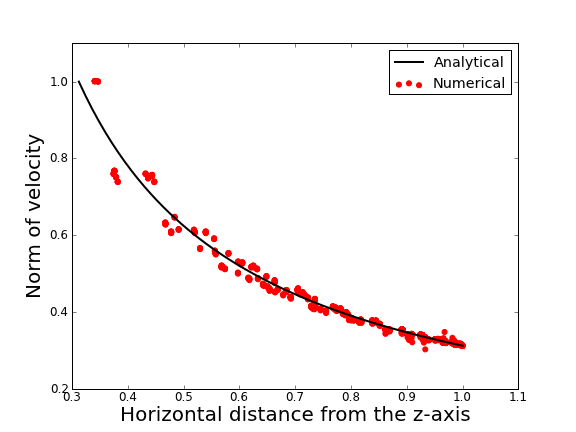}  & 
    \hspace{-0.02\textwidth}
    \includegraphics[scale=0.45]{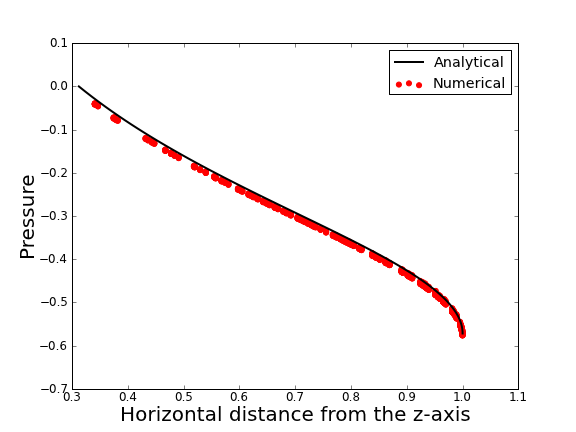} \\
    \hspace{-0.04\textwidth} 
    \includegraphics[scale=0.45]{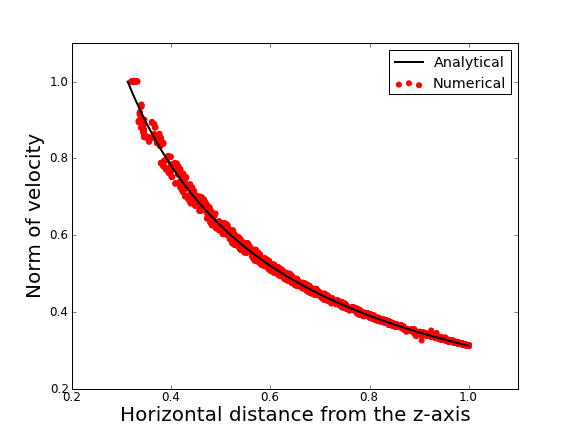}  & 
    \hspace{-0.02\textwidth}
    \includegraphics[scale=0.45]{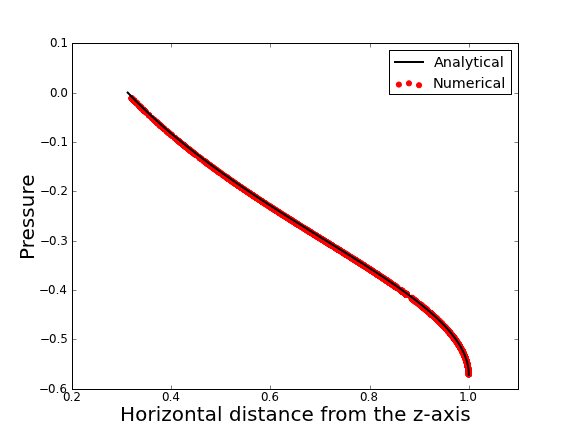}
  \end{tabular}
  \caption{A qualitative indication of convergence of the Whitney
    solution for the annular hemisphere on a series of
    non well-centered meshes similar to the mesh shown in the top
    panel of Figure~\ref{fig:hmsphrhlmsh}. Similar to
    Figure~\ref{fig:hmsphrhlFEEC}, the left panels are speed plots and
    the right ones are pressure plots. The rows of figures from top to
    bottom are profiles on increasingly finer meshes with 240, 1096
    and 4928 elements, respectively.}
  \label{fig:hmsphrhlFEECcnvrgnc}
\end{figure}

\begin{figure}[p]
 \begin{tabular}{m{0.475\textwidth}m{0.475\textwidth}}
    \hspace{-0.04\textwidth} \includegraphics[scale=0.45]
    {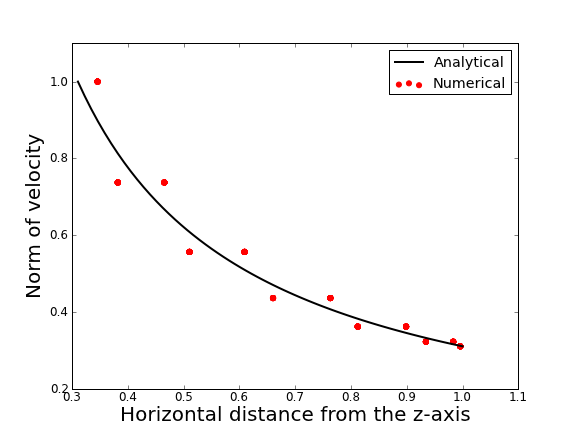}  & \hspace{-0.02\textwidth}
    \includegraphics[scale=0.45]{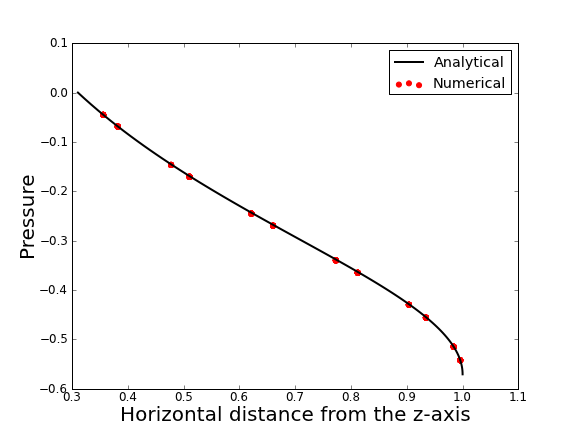}
    \\
    \hspace{-0.04\textwidth} \includegraphics[scale=0.45]
    {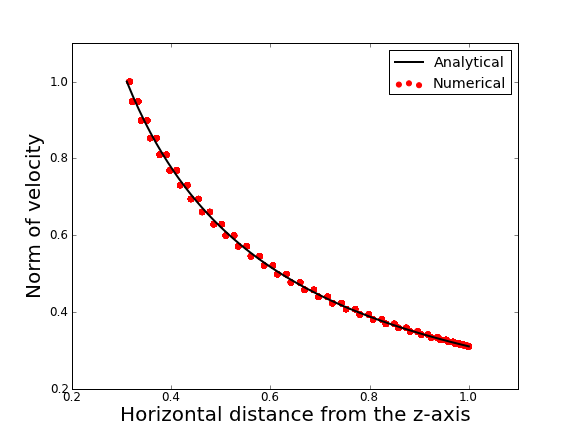}  & \hspace{-0.02\textwidth}
    \includegraphics[scale=0.45]{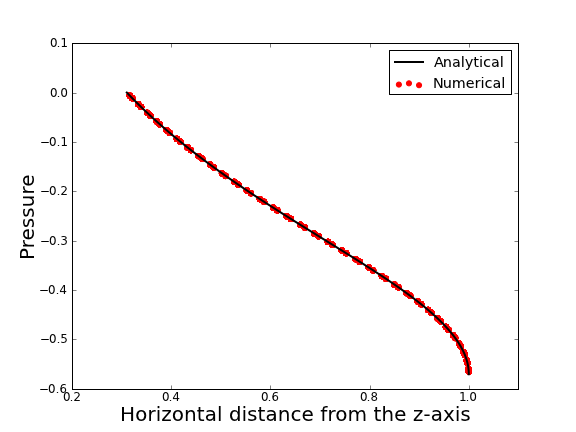}
    \\
    \hspace{-0.04\textwidth} \includegraphics[scale=0.45]
    {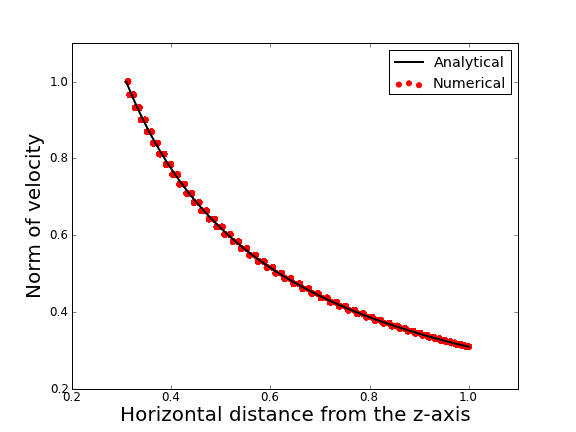}  & \hspace{-0.02\textwidth}
    \includegraphics[scale=0.45]{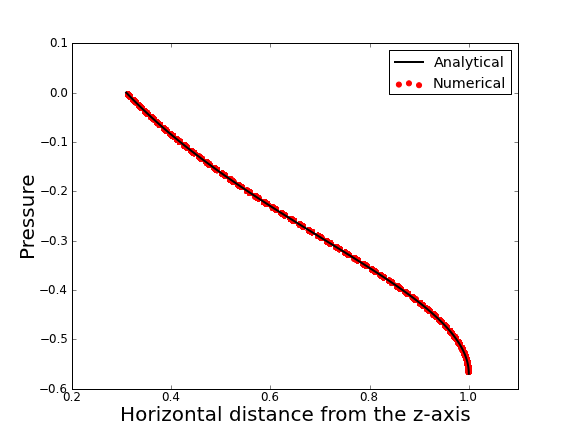}
  \end{tabular}
  \caption{A qualitative indication of convergence of DEC solution for
    the annular hemisphere on series of well-centered meshes similar
    to the one in the bottom panel of
    Figure~\ref{fig:hmsphrhlmsh}. Similar to
    Figure~\ref{fig:hmsphrhlFEEC}, the left panels are speed plots and
    the right ones are pressure plots. The rows of figures from top to
    bottom are profiles on increasingly finer meshes with 240, 6600
    and 15000 elements, respectively.}
  \label{fig:hmsphrhlDECcnvrgnc}
\end{figure}

\begin{figure}[p]
 \begin{tabular}{m{0.475\textwidth}m{0.475\textwidth}}
    \hspace{-0.04\textwidth} \includegraphics[scale=0.45]
    {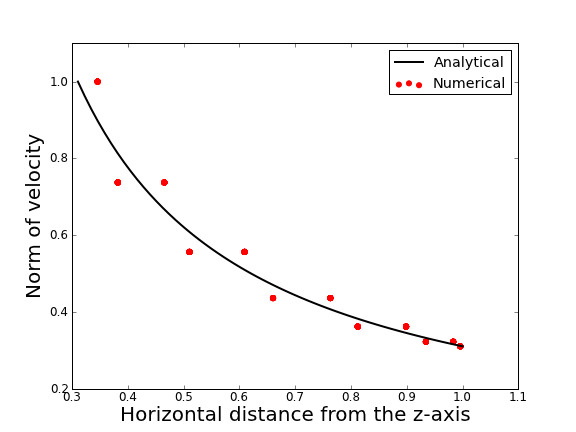}  & \hspace{-0.02\textwidth}
    \includegraphics[scale=0.45]{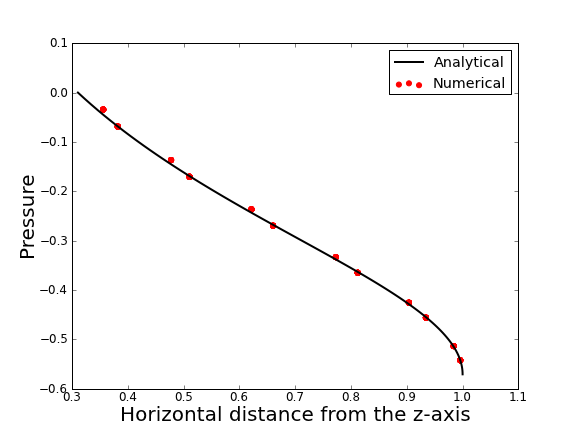}
    \\
    \hspace{-0.04\textwidth} \includegraphics[scale=0.45]
    {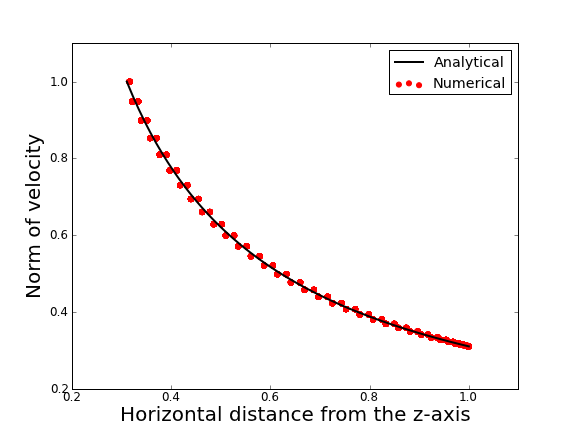}  & \hspace{-0.02\textwidth}
    \includegraphics[scale=0.45]{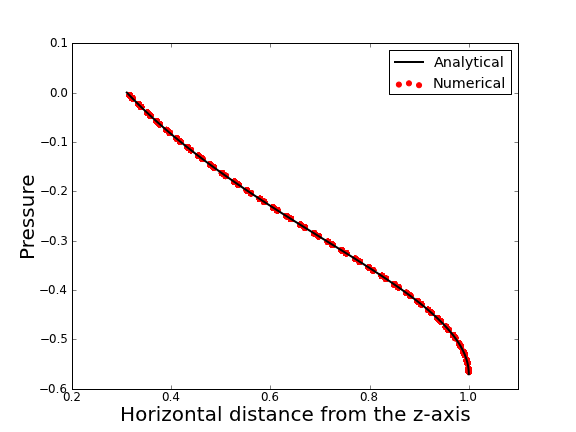}
    \\
    \hspace{-0.04\textwidth} \includegraphics[scale=0.45]
    {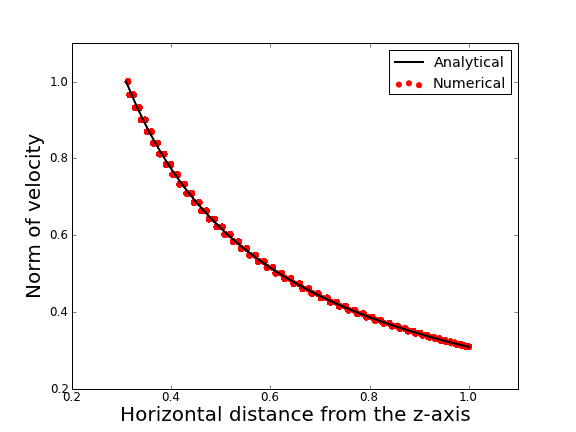}  & \hspace{-0.02\textwidth}
    \includegraphics[scale=0.45]{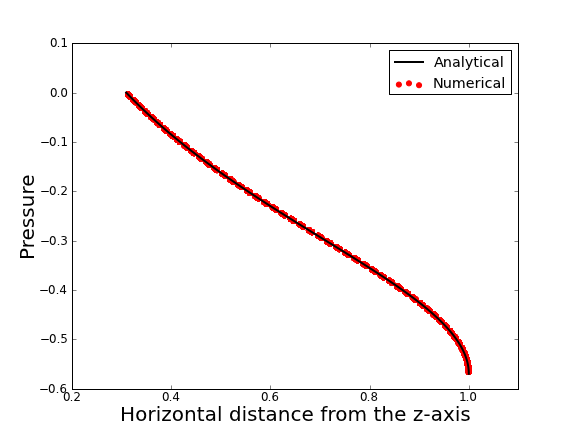}
  \end{tabular}
  \caption{A qualitative indication of convergence of the Whitney
    solution for the annular hemisphere on series of well-centered
    meshes similar to the one in the bottom panel of
    Figure~\ref{fig:hmsphrhlmsh}. Similar to
    Figure~\ref{fig:hmsphrhlFEEC}, the left panels are speed plots and
    the right ones are pressure plots. The rows of figures from top to
    bottom are profiles on increasingly finer meshes with 240, 6600
    and 15000 elements, respectively.}
  \label{fig:hmsphrhlFEECwctcnvrgnc}
\end{figure}

The results for the annular hemisphere experiments are shown in
Figures~\ref{fig:hmsphrhlFEEC}-\ref{fig:hmsphrhlFEECwctcnvrgnc}. For
every mesh used in these experiments the vertices always lie on the
hemisphere surface. For the DEC based method we use a well-centered
mesh, that is one in which each triangle is acute~\cite{VaHiGuRa2010a}
although a weaker condition can be imposed on the mesh. For the
Whitney solution method no such conditions are required. For
completeness we used well-centered as well as non well-centered meshes
when computing Whitney solutions. Some representative speed and
pressure plots for the Whitney solution are in
Figure~\ref{fig:hmsphrhlFEEC} and for DEC in
Figure~\ref{fig:hmsphrhlDEC}. For the purpose of visualizing the
computed pressures an $r$ value has to be associated with the computed
pressure corresponding to a triangle. For a well-centered mesh we use
the $r$ value of the circumcenter projected to the hemisphere and for
other meshes we use the projected barycenters. The flux is
interpolated and sampled at barycenter as described in the planar
case. For this flux, the $r$ value used is that of the barycenter
projected to the hemisphere.

\section{Conclusions and Future Work}

Our experience has been that it is easy to conduct numerical
experiments for Darcy flow on a surface with DEC and lowest order FEEC
using the Python package PyDEC~\cite{BeHi2011}. It would be
interesting to compare this with the programming effort required for a
vector calculus formulation on surfaces using classical mixed finite
elements. To continue the exterior calculus based experimentation, we
plan to do numerical experiments on the surface used here using
different boundary conditions, do computations using other surfaces,
and do quantitative convergence studies besides the qualitative
studies shown here.

\section*{Acknowledgment}
This work was funded in part by NSF CAREER Award Grant DMS-0645604. 

\bibliographystyle{acmdoi} 
\bibliography{hirani}

\begin{thebibliography}{1}

\bibitem{AbMaRa1988}
{\sc Abraham, R., Marsden, J.~E., and Ratiu, T.}
\newblock {\em Manifolds, Tensor Analysis, and Applications}, second~ed.
\newblock Springer--Verlag, New York, 1988.

\bibitem{ArFaWi2006}
{\sc Arnold, D.~N., Falk, R.~S., and Winther, R.}
\newblock Finite element exterior calculus, homological techniques, and
  applications.
\newblock In {\em Acta Numerica}, A.~Iserles, Ed., vol.~15. Cambridge
  University Press, 2006, pp.~1--155.

\bibitem{ArFaWi2010}
{\sc Arnold, D.~N., Falk, R.~S., and Winther, R.}
\newblock Finite element exterior calculus: from {H}odge theory to numerical
  stability.
\newblock {\em Bull. Amer. Math. Soc. (N.S.) 47}, 2 (2010), 281--354.
\newblock \href {http://dx.doi.org/10.1090/S0273-0979-10-01278-4} {\path{doi:
  10.1090/S0273-0979-10-01278-4}}.

\bibitem{BeHi2011}
{\sc Bell, N., and Hirani, A.~N.}
\newblock {PyDEC}: Algorithms and software for {D}iscretization of {E}xterior
  {C}alculus, March 2011.
\newblock Available as e-print on arxiv.org.
\newblock \href {http://arxiv.org/abs/1103.3076v1} {\path{arXiv:1103.3076v1}}.

\bibitem{HiNaCh2008}
{\sc Hirani, A.~N., Nakshatrala, K.~B., and Chaudhry, J.~H.}
\newblock Numerical method for {D}arcy flow derived using {D}iscrete {E}xterior
  {C}alculus.
\newblock Tech. Rep. UIUCDCS-R-2008-2937, Department of Computer Science,
  University of Illinois at Urbana-Champaign, 2008.
\newblock Also available as e-print on arxiv.org.
\newblock \href {http://arxiv.org/abs/0810.3434v3} {\path{arXiv:0810.3434v3}}.

\bibitem{VaHiGuRa2010a}
{\sc VanderZee, E., Hirani, A.~N., Guoy, D., and Ramos, E.~A.}
\newblock Well-centered triangulation.
\newblock {\em SIAM Journal on Scientific Computing 31}, 6 (2010), 4497--4523.
\newblock Also available as e-print on arxiv.org.
\newblock \href {http://arxiv.org/abs/0802.2108v3} {\path{arXiv:0802.2108v3}},
  \href {http://dx.doi.org/10.1137/090748214} {\path{doi: 10.1137/090748214}}.

\end{thebibliography}

\end{document}